\documentclass[11pt,fullpage]{article}

\newcommand{\ol}{\setlength{\itemsep}{0pt.}\begin{enumerate}}
\newcommand{\eol}{\end{enumerate}\setlength{\itemsep}{-\parsep}}
\newcommand{\ignore}[1]{}
\title{\bf A proof of hyperbolic van der Waerden
  conjecture : the right generalization is the ultimate simplification}
 
\author{Leonid Gurvits \thanks{%
{\tt gurvits@lanl.gov}. Los Alamos National Laboratory, 
Los Alamos, NM. } 
}

\begin{document}


\begin{titlepage}

\maketitle

\begin{abstract}
Consider a homogeneous polynomial $p(z_1,...,z_n)$ of degree $n$ in $n$ complex variables .
Assume that this polynomial satisfies the  property : \\

$|p(z_1,...,z_n)| \geq \prod_{1 \leq i \leq n} Re(z_i)$ on the domain $\{(z_1,...,z_n) :  Re(z_i) \geq 0 , 1 \leq i \leq n \}$ . \\

We prove that $|\frac{\partial^n}{\partial z_1...\partial z_n} p | \geq \frac{n!}{n^n}$ . \\
Our proof is relatively short and self-contained (i.e. we only use basic
properties of hyperbolic polynomials ). \\
As the van der Waerden conjecture for permanents , proved by D.I. Falikman and G.P. Egorychev ,
as well Bapat's conjecture for mixed discriminants , proved by the author ,
are particular cases of this result. \\ 
We also prove so called "small rank" lower bound (in the permanents context it corresponds to
sparse doubly-stochastic matrices , i.e. with small number of non-zero entries in each column).
The later lower bound generalizes (with simpler proofs) recent results by A.Schrijver for $k$-regular bipartite graphs.\\
Some important algorithmic applications are presented in the last section . 

\end{abstract} 
\end{titlepage}
\newpage

 
\newtheorem{THEOREM}{Theorem}[section]
\newenvironment{theorem}{\begin{THEOREM} \hspace{-.85em} {\bf :} 
}%
                        {\end{THEOREM}}
\newtheorem{LEMMA}[THEOREM]{Lemma}
\newenvironment{lemma}{\begin{LEMMA} \hspace{-.85em} {\bf :} }%
                      {\end{LEMMA}}
\newtheorem{COROLLARY}[THEOREM]{Corollary}
\newenvironment{corollary}{\begin{COROLLARY} \hspace{-.85em} {\bf 
:} }%
                          {\end{COROLLARY}}
\newtheorem{PROPOSITION}[THEOREM]{Proposition}
\newenvironment{proposition}{\begin{PROPOSITION} \hspace{-.85em} 
{\bf :} }%
                            {\end{PROPOSITION}}
\newtheorem{DEFINITION}[THEOREM]{Definition}
\newenvironment{definition}{\begin{DEFINITION} \hspace{-.85em} {\bf 
:} \rm}%
                            {\end{DEFINITION}}
\newtheorem{EXAMPLE}[THEOREM]{Example}
\newenvironment{example}{\begin{EXAMPLE} \hspace{-.85em} {\bf :} 
\rm}%
                            {\end{EXAMPLE}}
\newtheorem{CONJECTURE}[THEOREM]{Conjecture}
\newenvironment{conjecture}{\begin{CONJECTURE} \hspace{-.85em} 
{\bf :} \rm}%
                            {\end{CONJECTURE}}
\newtheorem{PROBLEM}[THEOREM]{Problem}
\newenvironment{problem}{\begin{PROBLEM} \hspace{-.85em} {\bf :} 
\rm}%
                            {\end{PROBLEM}}
\newtheorem{QUESTION}[THEOREM]{Question}
\newenvironment{question}{\begin{QUESTION} \hspace{-.85em} {\bf :} 
\rm}%
                            {\end{QUESTION}}
\newtheorem{REMARK}[THEOREM]{Remark}
\newenvironment{remark}{\begin{REMARK} \hspace{-.85em} {\bf :} 
\rm}%
                            {\end{REMARK}}
\newtheorem{FACT}[THEOREM]{Fact}
\newenvironment{fact}{\begin{FACT} \hspace{-.85em} {\bf :} 
\rm}%
		            {\end{FACT}}

 
\newcommand{\thm}{\begin{theorem}}
\newcommand{\lem}{\begin{lemma}}
\newcommand{\pro}{\begin{proposition}}
\newcommand{\dfn}{\begin{definition}}
\newcommand{\rem}{\begin{remark}}
\newcommand{\xam}{\begin{example}}
\newcommand{\cnj}{\begin{conjecture}}
\newcommand{\prb}{\begin{problem}}
\newcommand{\que}{\begin{question}}
\newcommand{\cor}{\begin{corollary}}
\newcommand{\fac}{\begin{fact}}

\newcommand{\prf}{\noindent{\bf Proof:} }
\newcommand{\ethm}{\end{theorem}}
\newcommand{\elem}{\end{lemma}}
\newcommand{\epro}{\end{proposition}}
\newcommand{\edfn}{\bbox\end{definition}}
\newcommand{\erem}{\bbox\end{remark}}
\newcommand{\exam}{\bbox\end{example}}
\newcommand{\ecnj}{\bbox\end{conjecture}}
\newcommand{\eprb}{\bbox\end{problem}}
\newcommand{\eque}{\bbox\end{question}}
\newcommand{\ecor}{\end{corollary}}
\newcommand{\efac}{\end{fact}}
\newcommand{\eprf}{\bbox}
\newcommand{\beqn}{\begin{equation}}
\newcommand{\eeqn}{\end{equation}}
\newcommand{\wbox}{\mbox{$\sqcap$\llap{$\sqcup$}}}
\newcommand{\bbox}{\vrule height7pt width4pt depth1pt}
\newcommand{\qed}{\bbox}

\newcommand{\rarrow}{\rightarrow}
\newcommand{\larrow}{\leftarrow}
\newcommand{\grad}{\bigtriangledown}

\overfullrule=0pt
\def\setof#1{\lbrace #1 \rbrace}

\section{Hyperbolic polynomials}

The following concept  of hyperbolic polynomials was originated in the theory of partial differential equations \cite{gar}, \cite{horm} ,\cite{kry} .\\
A homogeneous polynomial $p(x), x \in R^m$ of degree $n$ in $m$ real varibles is called hyperbolic in the direction $e \in R^m$ 
(or $e$- hyperbolic) if for any $x \in R^m$ the polynomial 
 $p(x - \lambda e)$  in the one variable $\lambda$ has exactly $n$ real roots counting their multiplicities. We assume in this
paper that $p(e) > 0$ .
Denote an ordered vector of roots of $p(x - \lambda e)$ as 
$\lambda(x) = (\lambda_{1}(x) \geq \lambda_{2}(x) \geq ... \lambda_{n}(x)) $. It is well known that the product of roots is equal
to $\frac{p(x)}{p(e)}$. Call $x \in R^m$ $e$-positive ($e$-nonnegative) if $\lambda_{n}(x) > 0$ ($\lambda_{n}(x) \geq 0$).
The fundamental result \cite{gar} in the theory of hyperbolic polynomials states that the set of $e$-nonnegative vectors is a 
closed convex cone. A $k$-tuple of vectors $(x_1,...x_k)$ is called $e$-positive ($e$-nonnegative) if $x_i , 1 \leq i \leq k$ are
$e$-positive ($e$-nonnegative). 
We denote the closed convex cone of $e$-nonnegative vectors as $N_{e}(p)$,
and the open convex cone of $e$-positive vectors
as $C_{e}(p)$.\\

{ \it Recent interest in the hyperbolic polynomials got sparked by the discovery \cite{gul} ,\cite{lewis} that $\log(p(x))$ is a self-concordant
barrier for the opened convex cone $C_{e}(p)$ and therefore the powerful mashinery of interior-point methods can be applied .
It is an important open problem whether this cone $C_{e}(p)$ has a semi-definite representation .}\\

It has been shown in \cite{gar} (see also \cite{khov}) that 
an $e$-hyperbolic polynomial $p$ is also
$d$- hyperbolic for all $e$-positive vectors $d \in C_{e}(p)$ ;
for all $d \in C_{e}(p)$ the set equalities $C_{d}(p) = C_{e}(p) , N_{d}(p) = N_{e}(p) $ . \\

Let us fix $n$ real vectors 
$x_i \in N_{e}(p) \subset R^m , 1 \leq i \leq n$ such that $\sum_{1 \leq i \leq n} x_i \in C_{e}(p)$
and define the following homogeneous polynomial:
\beqn
P_{x_1,..,x_n}(\alpha_1,...,\alpha_n) = p(\sum_{1 \leq i \leq n} \alpha_i x_i)
\eeqn
We will call such polynomials $P$-hyperbolic .\\
In other words ,a homogeneous polynomial $p(\alpha) , \alpha \in R^n$ of degree $n$ in $n$ real variables is $P$-hyperbolic
if it is $(1,1,..,1)$-hyperbolic ($e =(1,1,..,1)$)  and its closed cone of $e$-nonnegative vectors contains the nonnegative orthant 
$R^{n}_{+} = 
\{(x_1,...,x_n) : x_i \geq 0 , 1 \leq i \leq n \}$ . It follows from \cite{khov} that the coefficients of$P$-hyperbolic 
polynomials are nonnegative real numbers .\\

Following \cite{khov} , we define the $p$-mixed form of an $n$-vector tuple ${\bf X} = (x_1,..,x_n)$ as
\beqn
M_{p}({\bf X}) = : M_{p}(x_1,..,x_n) = \frac{\partial^n}{\partial \alpha_1...\partial \alpha_n} p(\sum_{1 \leq i \leq n} \alpha_i x_i)
\eeqn
Equivalently, the $p$-mixed form $M_{p}(x_1,..,x_n)$ can be defined by the polarization (see \cite{khov}) :
\beqn
M_{p}(x_1,..,x_n) = 2^{-n} \sum_{b_{i} \in \{-1, +1 \}, 1 \leq i \leq n }  p(\sum_{1 \leq i \leq n } b_i x_i) \prod_{1 \leq i \leq n } b_i
\eeqn

Associate with any vector $r = (r_1,...,r_n) \in I_{n,n}$ an $n$-tuple of $m$-dimensional vectors  ${\bf X}_{r}$ consisting
of  $r_i$ copies of $x_{i}  (1 \leq i \leq n) $.  
It follows from the Taylor's formula  that
\beqn
P_{x_1,..,x_n}(\alpha_1,...,\alpha_n) = \sum_{r \in I_{n,n}} \prod_{1 \leq i \leq n } \alpha_{i}^{r_{i}} M_{p}({\bf X}_{r})
\frac{1}{\prod_{1 \leq i \leq n } r_{i}!}
\eeqn
For an $e$-nonnegative tuple ${\bf X} = (x_1,..,x_n)$, define its capacity as:
\beqn
Cap({\bf X}) = \inf_{\alpha_i > 0, \prod_{1 \leq i \leq n }\alpha _i = 1} P_{x_1,..,x_n}(\alpha_1,...,\alpha_n)
\eeqn

Probably the best known example of a hyperbolic polynomial comes from the hyperbolic geometry :
\beqn
P(\alpha_0 ,...,\alpha_k) =\alpha_0^{2} - \sum_{1 \leq i \leq k } \alpha _i^{2}  
\eeqn

This polynomial is hyperbolic in the direction $(1,0,0,...,0)$.  Another "popular" hyperbolic polynomial is
$\det(X)$ restricted on a linear real space of hermitian $n \times n$ matrices .
In this case mixed forms are just mixed discriminants , hyperbolic direction is the identity matrix $I$ ,
the corresponding closed convex cone of $I$-nonnegative vectors coincides with a closed convex cone of positive semidefinite matrices .\\ 
Less known , but very interesting , hyperbolic polynomial is the Moore determinant $M\det(Y)$ 
restricted on a linear real space of hermitian quaternionic $n \times n$ matrices .
The Moore determinant is , essentially , the Pfaffian (see the corresponding definitions and
the theory in a very readable paper \cite{quat} ) .
The following definition is from \cite{half}.
\dfn
A polynomial $P(z_1,...,z_n)$ in $n$ complex variables is said to have the "half-plane property"
if $P(z_1,...,z_n) \neq 0$ provided $Re(z_i) > 0$ .
\edfn
In a control theory literature (see \cite{khar} ) the same property is called {\it Wide sense stability} .
And {\it Strict sense stability} means that \\
$P(z_1,...,z_n) \neq 0$ provided $Re(z_i) \geq 0$ .

The following simple fact shows that for homogeneous polynomials the "half-plane property" is , up to a
single factor , the same as $P$-hyperbolicity .
\pro
A homogeneous polynomial $R(z_1,...,z_n)$ has the "half-plane" property if and only if
the exists real $\alpha$ such that the polynomial $e^{i\alpha} R(z_1,...,z_n)$ is $P$-hyperbolic polynomial
with real nonnegative coefficients  .
\epro
\prf
\begin{enumerate}
\item
Suppose that $R(z_1,...,z_n) = e^{-i\alpha}Q(z_1,...,z_n)$ where $\alpha$ is real and $Q$ is $P$-hyperbolic.
Then $Q$ is $(1,1,...,)$-hyperbolic and all real vectors $(x_1,...,x_n)$ with positive coordinates are $(1,1,...,)$-positive .
Therefore $Q$ is $ (x_1,...,x_n)$-hyperbolic for all  real vectors $(x_1,...,x_n) \in R^{n}_{++}$ with positive coordinates .It follows
that $|R(x_1 + i y_1,...,x_n + i y_n)| = |Q(x_1 + i y_1,...,x_n + i y_n)| = |Q(x_1,...,x_n) \prod_{1 \leq k \leq n} (1 + i \lambda_{k})|$ ,
where $(\lambda_{1},...,\lambda_{n})$ are real roots of the real vector $(y_1,...,y_n)$ in the direction $(x_1,...,x_n)$.\\
This gives the following inequality , which is  equivalent to  the "half-plane property" of $R$ :
\begin{eqnarray}
&|R(x_1 + i y_1,...,x_n + i y_n)| \geq |R(x_1,...,x_n)| = \\ \nonumber
&= |Q(x_1,...,x_n)| > 0 : \\ \nonumber
&(x_1,...,x_n) \in R^{n}_{++} ,(y_1,...,y_n) \in R^{n} 
\end{eqnarray}
\item
Suppose that $R(z_1,...,z_n)$ has the "half-plane property" and consider the roots of the following polynomial equation in one complex variable :
$P(x_1 - z ,x_2 -z,...,x_n - z) = 0$ , where $(x_1,...,x_n) \in R^{n}$ is a real vector , $z = x +iy \in C$. If the imaginery part $Im(z) =y$
is not zero then , using the homogeniuty , $R( i \frac{x - x_{1}}{y} + 1,...,i \frac{x - x_{n}}{y} + 1) = 0$ , which is impossible as
$R$ has the "half-plane property". Therefore all roots of $R(X -t e) =0$ are real for all real vectors $X \in R^{n}$ (here $e = (1,1,...,1)$).
In the same way all roots of $R(X -t e) =0$ are real positive numbers if $X \in R^{n}_{++}$ . It follows that if $X \in R^{n}$ then
$R(X) = R(e) \prod{1 \leq k \leq n} \lambda_{k}(X)$ ,
where $(\lambda_{1},...,\lambda_{n})$ are ( real ) roots of the equation $R(X -t e) =0$ . Thus the polynomial $(\frac{1}{R(e)}) R$ takes
real values on $R^n$ and therefore its coefficients are real . In other words , the polynomial $(\frac{1}{R(e)}) R$ is $P$-hyperbolic .
If $R(1,1,...,1) = e^{- i \alpha} |R(1,1,...,1)|$ then the polynomial $e^{i\alpha} R$ is also $P$-hyperbolic .\\
(Recall that the coefficients of any $P$-hyperbolic polynomial $p$ are nonnegative for they are $p$-mixed forms of $e$-nonnegative tuples ,
and $p$-mixed forms of $e$-nonnegative tuples are nonnegative if $p(e) > 0$ \cite{khov}.)
\end{enumerate}
\eprf

\cor
Let $p(x_1,...,x_n)$ be a homogeneous polynomial in $n$ variables and of degree $n$ . Assume
that $p(1,1,...,1) > 0$ . Then the property \\
"polynomial $p$ is $P$-hyperbolic
and its capacity $Cap(p) =\inf_{x_i > 0, \prod_{1 \leq i \leq n }x_i = 1} p(x_1,..,x_n) = C > 0$ " \\
is equivalent to the property \\
"$\inf_{Re(z_i) > 0 ,\prod_{1 \leq i \leq n }Re(z_i) = 1} |p(z_1,...,z_n)| = C > 0$".
\ecor

\prf
Suppose that \\
 "polynomial $p$ is $P$-hyperbolic
and its capacity $Cap(p) =\inf_{x_i > 0, \prod_{1 \leq i \leq n }x_i = 1} p(x_1,..,x_n) = C > 0$ " .\\
Then , as in (7) , 
$$
|p(x_1 + i y_1,...,x_n + i y_n)| \geq p(x_1,...,x_n) \geq C \prod_{1 \leq i \leq n } x_i ; x_i \geq 0 , y_i \in R , 1 \leq i \leq n.
$$
Therefore $\inf_{Re(z_i) > 0 ,\prod_{1 \leq i \leq n }Re(z_i) = 1} |p(z_1,...,z_n)| = C$. \\

Assume that \\
"$\inf_{Re(z_i) > 0 ,\prod_{1 \leq i \leq n }Re(z_i) = 1} |p(z_1,...,z_n)| = C > 0$" . \\
Since $p(1,1,...,1) > 0$ , it follows from Proposition 1.2 that $p$ is $P$-hyperbolic . The equality \\
$Cap(p) = \inf_{Re(z_i) > 0 ,\prod_{1 \leq i \leq n }Re(z_i) = 1} |p(z_1,...,z_n)|$ follows .
\eprf

\rem
Corollary 1.3 essentially says that if $p(x_1,...,x_n)$ is a homogeneous polynomial in $n$ variables and of degree $n$
with real nonnegative coefficients and its complex capacity \\
$C-Cap(p) = :\inf_{Re(z_i) > 0 ,\prod_{1 \leq i \leq n }Re(z_i) = 1} |p(z_1,...,z_n)| = C >0$ \\
then its (real) capacity \\ 
$Cap(p) = \inf_{x_i > 0, \prod_{1 \leq i \leq n }x_i = 1} p(x_1,..,x_n) = C-Cap(p)$.\\
If $C-Cap(p) = 0$ then this statement can be wrong . I.e.
consider $q(x_1,...,x_n) = \frac{\sum_{1 \leq i \leq n} x_i^{n}}{n}$. Then $Cap(q)=1$ and $C-Cap(q) = 0$.
\erem

We use in this paper the following class of hyperbolic in the direction $(1,1,...,1)$  polynomials of degree $k$ :\\ 
$Q(\alpha_1,...,\alpha_k) = M_{p}(\sum_{1 \leq i \leq k} \alpha_i x_i,...,\sum_{1 \leq i \leq k} \alpha_i x_i , x_{k+1},...,x_{n}) $,
where $p$ is a  $e$-hyperbolic polynomial of degree $n > k$ , $(x_1,..,x_n)$ is $e$-nonnegative tuple , and the $p$-mixed form \\
$M_{p}(\sum_{1 \leq i \leq k} x_i,...,\sum_{1 \leq i \leq k} x_i , x_{k+1},...,x_{n}) > 0$. \\ 

\section{Main Theorem}
 \thm
\begin{enumerate}
\item
Let $q(x_1,x_2,...,x_n)$ be a $P$-hyperbolic (homogeneous) polynomial of degree $n$ .
Then
\beqn
\frac{\partial^n}{\partial x_1...\partial x_n} q(0,...,0) \geq  \frac{n!}{n^{n}} Cap(q)
\eeqn
\item
This bound is attained only on the following class of polynomials :
$$
q_{a_1,...,a_n}(x_1,...,x_n) =(\frac{\sum_{1 \leq i \leq n} a_i x_i }{n})^{n} ; a_i > 0 , 1 \leq i \leq n .
$$
(Notice that $Cap(q_{a_1,...,a_n}) = \prod_{1 \leq i \leq n} a_i$ .)

\end{enumerate}

\ethm
\subsection{Auxiliary Results}

\pro
\begin{enumerate}
\item
Let $c_1,...,c_n$ be real numbers ; $0 \leq c_i \leq 1 , 1 \leq i \leq n$ and $\sum_{1 \leq i \leq n} c_i = n-1$.\\
Define the following symmetric functions :
$$
S_{n}  =\prod_{1 \leq i \leq n} c_i , S_{n-1} = \sum_{1 \leq i \leq n} \prod_{j \neq i} c_j .
$$
Then the following entropic inequality holds :
$$
S_{n-1} - n S_{n} \geq e^{\sum_{1 \leq i \leq n} c_i \log(c_{i})} .
$$

\item (Mini van der Waerden conjecture)\\
Consider a doubly-stochastic $n \times n$ matrix $A= [a|b|...|b]$ .
I.e. $A$ has $n-1$ columns equal to the column vector $b$ , and one column equal to the column vector $a$ .
Let $a = (a_1,...,a_n)^{T} : a_i \geq 0 , sum_{1 \leq i \leq n} a_i = 1$ ;
$b = (b_1,...,b_n)^{T} : b_i = \frac{1 - a_i}{n-1} ,1 \leq i \leq n .$
Then the permanent $Per(A) \geq \frac{n!}{n^{n}}$ .
\end{enumerate}
\epro

\prf
\begin{enumerate}
\item
Doing simple "algebra" we get that
$$
S_{n-1} - n S_{n} = \prod_{1 \leq i \leq n} c_i (\sum_{1 \leq i \leq n} \frac{1-c_{i}}{c_i} ).
$$
Notice that $0 \leq 1-c_{i} \leq 1$ and $\sum_{1 \leq i \leq n} (1-c_i) = 1$.
Using the concavity of the logarithm we get that
$$
\log(S_{n-1} - n S_{n}) \geq \sum_{1 \leq i \leq n} \log(c_i) + \sum_{1 \leq i \leq n} (1-c_{i}) \log(\frac{1}{c_i}) = 
\sum_{1 \leq i \leq n} c_i \log(c_{i}).
$$
\item
$$
per(A) = \frac{ (n-1)!}{ (n-1)^{n-1}} \sum_{1 \leq i \leq n} a_i \prod_{j \neq i} (1-a_j) .
$$
Define $c_i = 1-a_i$ . Then $0 \leq 1-c_{i} \leq 1$ , $\sum_{1 \leq i \leq n} c_i = n-1$ and the permanent
$$
Per(A) = \frac{ (n-1)!}{ (n-1)^{n-1}} (S_{n-1} - n S_{n} ) .
$$
It is easy to prove and well known that
$$
\min_{0 \leq 1-c_{i} \leq 1 ;\sum_{1 \leq i \leq n} c_i = n-1 } \sum_{1 \leq i \leq n} c_i \log(c_{i})  = \sum_{1 \leq i \leq n} 
\frac{n-1}{n} \log (\frac{n-1}{n}) =
\log ((\frac{n-1}{n})^{n-1}) .
$$
Using the entropic inequality from the first part we get the following equality
$$
\min_{0 \leq 1-c_{i} \leq 1 ;\sum_{1 \leq i \leq n} c_i = n-1 } S_{n-1} - nS_{n} =
(\frac{n-1}{n})^{n-1} .
$$

Which gives the needed inequality
$$
Per(A) \geq \frac{ (n-1)!}{ (n-1)^{n-1}} (\frac{n-1}{n})^{n-1} = \frac{n!}{n^{n}} .
$$

It is easy to see (strict concavity of $\sum_{1 \leq i \leq n} c_i \log(c_{i})$)
that the last inequality is strict unless $ A(i,j) = \frac{1}{n} ; 1 \leq i,j \leq n$.
\end{enumerate}
\eprf

\cor
Define capacity of $n \times n$ matrix $A$ with nonnegative entries
as 
$$
Cap(A) = \inf_{x_j > 0} \frac{\prod_{1 \leq i \leq n} \sum_{1 \leq j \leq n} A(i,j)x_j} {\prod_{1 \leq j \leq n} x_j}
$$
If $A = [c|d|...|d]$ then $Per(A) \geq \frac{n!}{n^{n}} Cap(A)$ .

\ecor

\prf
Sinkhorn's diagonal scaling to doubly-stochastic matrices does the job. I.e. ,
if $A = [c|d|...|d]$  and all entries of $A$ are positive then
there exist two diagonal matrices with positive entries $D_1 , D_2$
and a doubly-stochastic matrix $B= [a|b|...|b]$
such that $\det(D_1 D_2) = Cap(A)$ and $A = D_1 B D_2$ .
\eprf

\cor
Consider an univariate polynomial \\
$R(t) = \sum_{0 \leq i \leq n} d_i t^{i} = \prod_{1 \leq i \leq n} (a_i t + b_i)$ ,
where $a_i , b_i \geq 0$ . If for some positive real number $C$ the inequality 
$R(t) \geq C t$ holds for all $t \geq 0$ then
\beqn
d_1 = \frac{\partial}{\partial t} R(0) \geq C((\frac{n-1}{n})^{n-1})
\eeqn
The inequality (9) is attained on the polynomial $R(t) = n^{-n} (t + n- 1)^{n} .$
\ecor

\prf
Associate with polynomial $R(t) = \prod_{1 \leq i \leq n} (a_i t + b_i)$ the following
matrix $A = [a |c|...|c] $ , where $a = (a_1,...,a_n)^{T} , c=\frac{1}{n-1} (b_1,...,b_n)^{T}$.
The condition $R(t) \geq C t , \forall t \geq 0$ is equivalent to the inequality
$Cap(A) \geq C$ . And $d_1 = \frac{(n-1)!} {(n-1)^{n-1}} Per(A) .$ It follows from Corollary 2.3 that
$$
d_1 = (\frac{(n-1)!} {(n-1)^{n-1}} )^{-1} Per(A) \geq (\frac{(n-1)!} {(n-1)^{n-1}})^{-1} (\frac{n!}{n^{n}} C) = C((\frac{n-1}{n})^{n-1}).
$$
\eprf

\pro
Let $p(X)$ be $e$-hyperbolic (homogeneous) polynomial of degree $n$ , $p(e) > 0$ . 
Consider two $e$-nonnegative vectors
$Z,Y \in N_{e}(p)$ such that $Z+Y \in C_{e}(p)$ , i.e. $Z+Y$ is $e$-positive .
Then
\beqn
p(tZ +Y) = \prod_{1 \leq i \leq n} (a_i t + b_i ) ; a_i,b_i \geq 0 , a_i + b_i > 0 , 1 \leq i \leq n .
\eeqn
\epro

\prf
As the vector $Z+Y = D$ is $e$-positive hence $p(Z+Y) > 0$ ,the polynomial $p$ is $Z+Y$-hyperbolic
and any $e$-positive ($e$-nonnegative) is also $Z+Y$-positive($Z+Y$-nonnegative) \cite{khov} .
Doing simple algebra , we get that $p(tZ +Y) = p((t-1)Z +D)$. \\

Let $0 \leq \lambda_{1} \leq \lambda_{2} \leq ...\leq \lambda_{n} $ be nonnegative roots
of the equation $p(Z- x D)) = 0$. Since $D-X = Y \in N_{e}(p) = N_{D}(p) $ hence $\lambda_{n} \leq 1$.
Therefore
$$
p(tZ +Y) = p((t-1)Z +D) = p(D) \prod_{1 \leq i \leq n} (t \lambda_{i} + (1-\lambda_{i}) )
$$
We can put $a_i = (p(Z+Y)) \lambda_{i}  \geq 0 , b_i = (p(Z+Y)) (1-\lambda_{i}) \geq 0$.
\eprf

\pro
Let $q(x_1,x_2,...,x_n)$ be a $P$-hyperbolic (homogeneous) polynomial of degree $n$ .
Define a new homogeneous polynomial of degree $n-1$ in $n-1$ variables :
$$
r(x_2,...,x_n) = \frac{\partial}{\partial x_1} q(0,x_2,...,x_n) .
$$ 
If $Cap(q) > 0$ then the polynomial $r$ is also $P$-hyperbolic .
\epro

\prf
Proved in \cite{newhyp} , easy modification of the argument in \cite{khov} , essentially
the Rolle's theorem .
\eprf

The next Lemma is the final auxiliary Result .

\lem

Define $F(n) = \frac{n!}{n^n}$ . The following inequality holds :
\beqn
Cap(r) \geq \frac{F(n)}{F(n-1)} Cap(q) = ((\frac{n-1}{n})^{n-1} )Cap(q)
\eeqn
\elem

\prf
Fix positive real numbers $(x_2,...,x_n)$ such that $\prod_{2 \leq i \leq n} x_i = 1.$
Define the following two real $n$-dimensionals vectors with nonnegative coordinates :
$Z= (1,0,0,...,0) , Y = (0,x_2,...,x_n)$ . The vector $Z+Y$ is $e$-positive .
Consider the next univariate polynomial $R(t) = p(tZ +Y)$.
It follows from Proposition 2.5 that
$$
R(t) = \prod_{1 \leq i \leq n} (a_i t + b_i) = \sum_{0 \leq i \leq n} d_i t^{i} ,
$$
where $a_i , b_i \geq 0$  and $r(x_2,...,x_n)= d_1$. \\
We get from the definition of $Cap(q)$ that 
$$
R(t) = \prod_{1 \leq i \leq n} (a_i t + b_i) = p(t,x_2,...,x_n) \geq Cap(q) t\prod_{2 \leq i \leq n} x_i  = t Cap(q) .
$$
Using Corollary 2.4 , we get that 
$$
r(x_2,...,x_n)= d_1 \geq ((\frac{n-1}{n})^{n-1} )Cap(q) .
$$
In other words , that $Cap(r) \geq \frac{F(n)}{F(n-1)} Cap(q) = (\frac{n-1}{n})^{n-1} )Cap(q)$.

\subsection{Proof of the Main Theorem}

(Only first part of Theorem 2.1 is proved in this draft .
The uniqueness part will be presented in the final version .)\\ 
\prf
Our proof is by (simple and natural) induction in $n$ . Theorem 2.1 is clearly true for
$n=1$ . Suppose it is true for all $k \leq n-1$ .
Let $q(x_1,x_2,...,x_n)$ be a $P$-hyperbolic (homogeneous) polynomial of degree $n$ and $Cap(q) = C > 0$.
Then using Lemma 2.7 we get that 
$$
Cap(r) \geq((\frac{n-1}{n})^{n-1} ) C = \frac{F(n)}{F(n-1)} C ,
$$
where $F(n) = \frac{n!}{n^n}$ and $r(x_2,...,x_n) = \frac{\partial}{\partial x_1} q(0,x_2,...,x_n)$ is 
a $P$-hyperbolic (homogeneous) polynomial of degree $n-1$ . Using induction
we get the needed inequality
$$
\frac{\partial^n}{\partial x_1...\partial x_n} q(x_1,...,x_n) = 
\frac{\partial^{n-1}}{\partial x_2...\partial x_n} r(x_2,...,x_n) \geq F(n-1) Cap(r) \geq F(n-1)\frac{F(n)}{F(n-1)} Cap(q) =
\frac{n!}{n^n} Cap(q) .
$$
\eprf

\xam
Consider a $n$-tuple of quaternionic hermitian $n \times n$ matrices ${\bf H} =(H_1,...,H_n)$
and define the following homogeneous polynomial of degree $n$ in $n$ real variables :
$$
Q_{{\bf H}}(x_1,...,x_n) = M\det(\sum_{1 \leq i \leq n} x_i H_i) ,
$$
where $M\det$ is the Moore determinant (consult the fantastic survey \cite{quat} 
on the subject of various quaternionic determinants). 
It is well known that right eigenvalues of quaternionic hermitian matrices are real (in this case
the Moore's determinant is equal to the product of right eigenvalues ) ,
quaternionic hermitian matrices with all right eigenvalues being nonnegative called
quaternionic positive semidefinite ( we write $H \succeq 0$ if the quaternionic hermitian matrix
$H$ is quaternionic positive semidefinite .) The
$tr(H)$ is equal to the sum of all (real) right eigenvalues of $H$ .\\
A $n$-tuple of quaternionic hermitian $n \times n$ matrices ${\bf H} =(H_1,...,H_n)$ is called
doubly stochastic if :
$$
H_i \succeq 0 , tr (H_i) = 1 , 1 \leq i \leq n ; \sum_{1 \leq i \leq n}  H_i = I .
$$
It is straigthforward to prove that if the tuple ${\bf H} =(H_1,...,H_n)$ is doubly stochastic
then the polynomial $Q_{{\bf H}}(x_1,...,x_n)$ is $P$-hyperbolic and $Cap(Q_{{\bf H}}) =1$.
It follows  from Theorem 2.1 that if the tuple ${\bf H} =(H_1,...,H_n)$ is doubly stochastic
then the following inequlity holds :
\beqn
HM({\bf H}) =:\frac{\partial^n}{\partial x_1...\partial x_n} M\det(\sum_{1 \leq i \leq n} x_i H_i) \geq \frac{n!}{n^n}
\eeqn 

If the tuple ${\bf H} =(H_1,...,H_n)$ consists of real diagonal positive semidefinite matrices 
then inequality (12) is the statement of the van der Waerden conjecture for permanents proved in \cite{fal} ;
if the tuple ${\bf H} =(H_1,...,H_n)$ consists of complex hermitian positive semidefinite matrices
then inequality (12) is the statement of the Bapat's conjecture \cite{bapat} for mixed discriminants
proved by the author in \cite{gur} . Even this quaternionic case seems to be a new result .
\exam

\rem
Notice that we did not use Falikman-Egorychev theorem (\cite{fal} , \cite{ego}) which proves
the "first" van der Waerden Conjecture \cite{minc} , bur rather its particularly
simple case (Proposition 2.2) . Theorem 2.2 generalizes all known variants
of  van der Waerden Conjecture (\cite{bapat} , \cite{gur} and others ...).
It also proves as Hall's theorem on perfect bipartite matchings , Rado's theorem
and its hyperbolic analogue \cite{newhyp} , \cite{hyp} .
And we did not use the Alexandrov-Fenchel inequalities ... \\
The main "spring" of our proof is that we work in a very large class
of $P$-hyperbolic polynomials , this class is large enough to allow the
easy induction. In fact , the clearest (in our opinion) proof of the
Alexandrov-Fenchel inequalities for mixed discriminants is
in A.G. Khovanskii' 1984 paper \cite{khov} . The Khovanskii' proof
is based on the similar induction (via partial differentions) to
the one used in our paper . In a way , the Alexandrov-Fenchel inequalities 
are "hidden" in our proof . \\

\erem

\subsection{Small Rank Lower Bound}
\dfn
Consider  a homogeneous polynomial $p(x) , x \in R^m$ of degree $n$ in $m$ real variables which is
hyperbolic in the direction $e$. Denote an ordered vector of roots of $p(x - \lambda e)$ as 
$\lambda(x) = (\lambda_{1}(x) \geq \lambda_{2}(x) \geq ... \lambda_{n}(x)) $ . We
define the $p$-rank of $x \in R^m$ in direction $e$ as $Rank_{p}(x) = |\{i : \lambda_{i}(x) \neq 0 \} |$.
It follows from Theorem 1.5 that the $p$-rank of $x \in R^m$ in any direction $d \in C_{e}$ is equal to the $p$-rank of $x \in R^m$
in direction $e$ , which we call the $p$-rank of $x \in R^m$ .
\edfn
Consider the following polynomial in one variable $D(t) = p(td+x) = \sum_{0 \leq i \leq n} c_{i} t ^{i} $.
It follows from the identity (4) that \\
\begin{eqnarray}
&c_{n} = M_{p}(d,..,d) (n!)^{-1} = p(d), \\ \nonumber
&c_{n-1} = M_{p}(x,d,..,d) (1! (n-1)!)^{-1},...,\\ \nonumber
&c_{0} = M_{p}(x,..,x) (n!)^{-1} = p(x).
\end{eqnarray}
Let $(\lambda_{1}^{(d)}(x) \geq \lambda_{2}^{(d)}(x) \geq... \geq \lambda_{n}^{(d)}(x))$ be the (real) roots of $x$
in the $e$-positive direction $d$, i.e. the roots of the equation $p(td-x) = 0$ .
Define (canonical symmetric functions) :
$$
S_{k,d}(x) = \sum_{1 \leq i_1 < i_2 <  ...< i_k \leq n} \lambda_{i_{1}}(x) \lambda_{i_{2}}(x)... \lambda_{i_{k}}(x).
$$
Then $S_{k,d}(x) = \frac{c_{n-k}}{c_{n}}$ . Clearly if $x$ is $e$-nonnegative then for any e-positive vector $d$
the $p$-rank $Rank_{p}(x) =\max \{k : S_{k,d}(x) > 0 \}$ .
The following usefull result can be found in \cite{khov} (the proof is essentially the same induction via partial differentions).

\fac 
Consider a homogeneous polynomial $p(x) , x \in R^m$ of degree $n$ in $m$ real variables which is
hyperbolic in the direction $e , p(e) > 0$. Then the following statements are true :
\begin{enumerate}
\item
The $p$-mixed form $M_{p}(y_1,...,y_n)$ is linear in each $y_i \in R^M$  when the rest is fixed .
\item
If the vectors $y_i , 1 \leq i \leq n$ are $e$-positive ($e$-nonnegative) then $M_{p}(y_1,...,y_n) > 0$ ($M_{p}(y_1,...,y_n) \geq 0$).
\item
If the vectors $y_i , z_i , y_i-z_i \in R^M: ; 1 \leq i \leq n$ are $e$-nonnegative then 
$$
M_{p}(y_1,...,y_n) \geq M_{p}(z_1,...,z_n).
$$
\end{enumerate}
\efac 

One of the corollaries of this fact is that for $e$-nonnegative vectors $x$ the number of positive roots of the univariate
equation $p(td-x) = 0$ is the same for all $e$-positive vectors $d$.
\pro
Let $q(x_1,x_2,...,x_n)$ be a $P$-hyperbolic (homogeneous) polynomial of degree $n$ and $Cap(q) > 0$.
Define a new homogeneous polynomial of degree $n-1$ in $n-1$ variables :
$$
r(x_2,...,x_n) = \frac{\partial}{\partial x_1} q(0,x_2,...,x_n) .
$$ 
Let $(e_1,e_2,...,e_n)$ be a canonical basis in $R^n$ . In other words , the vector $e_i \in R^n$ is the $i$th column
of $n \times n$ identity matrix $I$ . Then for all $2 \leq i \leq n$ the following inequality holds
\beqn
 Rank_{r}(e_{i}) \leq \min(Rank_{q}(e_{i}) , n-1)
\eeqn
\epro

\prf
First we recall the following formula , expressing the polynomial $\frac{\partial}{\partial x_1} q(0,x_2,...,x_n)$
in terms of $q$-mixed forms (\cite{khov},\cite{newhyp} ) :
$$
r(x_2,...,x_n) = M_{q}(e_1, z,...,z)  ((n-1)!)^{-1} , z = (0,x_2,...,x_n)^{T} .
$$
Clearly ,  $Rank_{r}(e_{i}) \leq n-1 \leq \min(Rank_{q}(e_{i}) , n-1)$ if $Rank_{q}(e_{i}) \geq n-1$.
Suppose that $Rank_{q}(e_{i}) = R_i \leq n-2$ . Since the vectors $(e_1,e_2,...,e_n)$ are $e$-nonnegative
hence
$$
M_{q}(e,...,e,e_i,...,e_i) =0 ,
$$
where the $n$-tuple $(e,...,e,e_i,...,e_i)$ contains $R_i +1 = Rank_{r}(e_{i})+1$ copies of $e_i$ and $n-1 -R_i$ copies
of $e = (1,1,...,1)^{T}$.
Define $d = \sum_{2 \leq i \leq n} e_i = e -e_1$.
To prove that $Rank_{r}(e_{i}) \leq \min(Rank_{q}(e_{i}) , n-1)$ we need to prove
that $M_{r}(d,..,d,e_i,...,e_i) = 0$ , where the $n-1$-tuple $(d,..,d,e_i,...,e_i)$ contains
$R_i +1 = Rank_{r}(e_{i})+1$ copies of $e_i$ and $n-2 -R_i$ copies
of $d$. But 
$$
M_{r}(d,..,d,e_i,...,e_i) = M_{q}(e_1,d,...,d ,e_i,...,e_i) .
$$
We have now two $n$-tuples $ {\bf T}_{1} = (e,...,e,e_i,...,e_i)$ and $ {\bf T}_{2} =(e_1,d,...,d ,e_i,...,e_i)$ .
The $n$-tuples  ${\bf T}_{1} ,{\bf T}_{2} , {\bf T}_{1} - {\bf T}_{2}$ consist of
$e$-nonnegative vectors . Therefore , using the monotonicity result from \cite{khov} ,
we get that $M_{r}(d,..,d,e_i,...,e_i) \leq M_{q}(e,...,e,e_i,...,e_i) = 0$.
\eprf

\lem
Let $q(x_1,x_2,...,x_n)$ be a $P$-hyperbolic (homogeneous) polynomial of degree $n$ and
$Rank_{q} (e_{1}) = k$ .
Then following inequality holds :
\beqn
Cap(r) \geq  ((\frac{k-1}{k})^{k-1} )Cap(q)
\eeqn
\elem

\prf (Very similar to the proof of Lemma 2.7).\\

Fix positive real numbers $(x_2,...,x_n)$ such that $\prod_{2 \leq i \leq n} x_i = 1.$
Define the following two real $n$-dimensionals vectors with nonnegative coordinates :
$Z= (1,0,0,...,0) , Y = (0,x_2,...,x_n)$ . The vector $Z+Y$ is $e$-positive .
Consider the next univariate polynomial $R(t) = p(tZ +Y)$.
It follows from Proposition 2.5 that
$$
R(t) = \prod_{1 \leq i \leq n} (a_i t + b_i) = \sum_{0 \leq i \leq n} d_i t^{i} ,
$$
where $a_i , b_i \geq 0$  and $r(x_2,...,x_n)= d_1$ and
the cardinality $|\{i : a_i > 0 \}| = k$.  In other words the degree
$deg(R) = k$\\
We get from the definition of $Cap(q)$ that 
$$
R(t) = \prod_{1 \leq i \leq n} (a_i t + b_i) = p(t,x_2,...,x_n) \geq Cap(q) t\prod_{2 \leq i \leq n} x_i  = t Cap(q) .
$$
Using Corollary 2.4 , we get that 
$$
r(x_2,...,x_n)= d_1 \geq ((\frac{k-1}{k})^{k-1} )Cap(q) .
$$
In other words , that $Cap(r) \geq (\frac{k-1}{k})^{k-1} )Cap(q)$.
\eprf

\thm
\begin{enumerate}
\item
Let $q(x_1,x_2,...,x_n)$ be a $P$-hyperbolic (homogeneous) polynomial of degree $n$ ;
$Rank_{q} (e_{i}) = R_i$ . Define $G_i = \min(R_{i} , n+1-i)$
Then
\beqn
\frac{\partial^n}{\partial x_1...\partial x_n} q(0,...,0) \geq  \prod_{1 \leq i \leq n}  (\frac{G_{i} -1}{G_{i}})^{G_{i} -1} Cap(q)
\eeqn
\item
If $Rank_{q} (e_{i}) \leq k \leq n$ then
\beqn
\frac{\partial^n}{\partial x_1...\partial x_n} q(0,...,0) \geq (\frac{k -1}{k})^{(k -1)(n-k)} \frac{k!}{k^{k}} Cap(q)
\eeqn
\end{enumerate}

\ethm

\prf
We use the same induction as in the proof of Theorem 2.1 together with Proposition 2.12 and Lemma 2.13 .
\eprf

The following result is a direct corollary Theorem 2.14 . Even the permanental inequality (18) seems
to be new (compare (18) with the corresponding result from \cite{schr} ). The easiness of our proof
(compare again with \cite{schr} ) suggests that the "method of hyperbolic polynomials" introduced in this paper is
very powerful and natural .
\cor
\begin{enumerate}
\item
Consider a doubly-stochastic $n$-tuple ${\bf A} = (A_1,...,A_n)$ of $n \times n$ hermitian positive semidefinite matrices , \\
i.e. $A_i \succeq 0 , tr(A_i = 1 ; 1 \leq i \leq n$ and $\sum_{1 \leq i \leq n} A_i = I$ .\\
If $Rank(A_i) \leq k \leq n$ then the mixed discriminant
$$
M(A_1,...,A_n) = :\frac{\partial^n}{\partial x_1...\partial x_n} \det(\sum_{1 \leq i \leq n} x_i A_i)  \geq (\frac{k -1}{k})^{(k -1)(n-k)} \frac{k!}{k^{k}} .
$$
\item
Let $A = \{A(i,j) : 1 \leq i,j \leq n\}$ be a doubly-stochastic $n \times n$ matrix . Suppose that
the cardinalities $| \{j : A(i,j) > 0\}| \leq k \leq n$ for $1 \leq i \leq n-k$ . Then the following
permanental inequality holds   :
\beqn
Per(A) \geq (\frac{k -1}{k})^{(k -1)(n-k)} \frac{k!}{k^{k}}.
\eeqn
\end{enumerate}
\ecor

\section{Applications}
Suppose that a $P$-hyperbolic (aka {\it Strict sense stable} homogeneous polynomial)
$$
p(x_1,...,x_n) =\sum_{ \sum_{1 \leq i \leq n} r_{i} =n } a_{(r_1,...,r_n) } \prod_{1 \leq i \leq n} x_{i}^{r_{i}}
$$ 
has nonnegative integer components coefficients and given as an oracle . I.e. we don't have a list coefficients ,
but can evaluate $p(x_1,...,x_n)$ on rational inputs . \\
An algorithm is called  deterministic polynomial-time oracle if
it evaluates the given polynomial $p(.)$ at a number of rational vectors $(q_1,...,q_n)$ which is 
polynomial in $n$ and $\log(p(1,1,..,1))$; these rational vectors $(q_1,...,q_n)$ are supposed to 
have bit-wise complexity which is polynomial in $n$ and $\log(p(1,1,..,1))$ ; and the additional auxilary arithmetic 
computations also take a polynomial number of steps in $n$ and $\log(p(1,1,..,1))$ .\\ 

The following theorem combines the algorithm from \cite{newhyp} and Theorem 2.1 .

\thm
There exists a deterministic polynomial-time oracle algorithm which computes for given as an oracle
$P$-hyperbolic polynomial $p(x_1,...,x_n)$ a number $F(p)$ satisfying the inequality
$$
\frac{\partial^n}{\partial x_1...\partial x_n} p(0,...,0) \leq F(p) \leq 
e^{n} \frac{\partial^n}{\partial x_1...\partial x_n} p(x_1,...,x_n) .
$$
\ethm

Theorem 3.1 can be (slightly) improved . I.e. it can be applied
to the polynomial 
$$
p_{k}(x_{k+1},...,x_{n}) =\frac{\partial^k}{\partial x_1...\partial x_k} p(0,..,0, x_{k+1},...,x_{n}).
$$
Notice that the polynomial $p_{k}$ is a homogeneous polynomial of degree $n-k$ in $n-k$ variables .
It is easy to prove that if $p= p_{0}$ is $P$-hyperbolic and $Cap(p) > 0$ then for all $0 \leq k \leq n$
the polynomials $p_{k}$ are also $P$-hyperbolic and $Cap(p_{k}) > 0$ . \\
The trick is that if $k = m \log_{2}(n)$ then
the polynomial $p_{k}$ can be evaluated using $O(n^{m+1})$ oracle calls of the (original) polynomial
$p$ . This observations allows to decrease the multiplicative factor in Theorem 3.1 from
$e^{n}$ to $\frac{e^{n}}{n^{m}}$ for any fixed $m$ . If the polynomial $p$ can
be explicitly evaluated in deterministic polynomial time , this observation
results in deterministic polynomial time algorithms to approximate
$\frac{\partial^n}{\partial x_1...\partial x_n} p(0,...,0)$ within multiplicative
factor $\frac{e^{n}}{n^{m}}$ for any fixed $m$ . Which is an improvement
of results in \cite{lsw} (permanents , $p$ is a multilinear polynomial) and in \cite{GS} , \cite{GS1} (mixed discriminants
$p$ is a determinantal polynomial) .

\section{Open Problems and Acknowledgements}

\prb
Is first part of Theorem 2.1 true for the volume polynomials
$p(x_1,...,x_n) = Volume(\sum_{1 \leq i \leq n} x_i C_i) $ ,
where $C_i ,1 \leq i \leq n$ are convex compact subsets of $R^{n}$ ? \\
Not all volume polynomials are $P$-hyperbolic (see the example in \cite{khov}) .
\eprb

\prb
What is a "good" model of a random $P$-hyperbolic polynomial ?
By "good" we mean that with high probability the inequality (8)
is much tighter . I.e. with high probability 
$$
\frac{\partial^n}{\partial x_1...\partial x_n} q(0,...,0) \geq  (1 + O(n^{-1})) Cap(q)
$$
\eprb

After the first draft had been posted Hugo Woerdeman found more direct proof of Corollary 2.4 .\\
I would like to thank Mihai Putinar , Sergey Fomin , George Soules , Alex Samorodnitsky  for
the interest to this paper .

\end{document}